\documentclass[12pt]{amsart}
\usepackage{graphicx}
\usepackage{amssymb}
\usepackage{amsfonts}
\usepackage{amsmath}
\usepackage{array}
\usepackage{rotating}

\usepackage[matrix,frame]{xypic}
\newdimen\vcadre\vcadre=0.1cm 
\newdimen\hcadre\hcadre=0.1cm 
\def\GrTeXBox#1{\vbox{\vskip\vcadre\hbox{\hskip\hcadre%
      $#1$%
   \hskip\hcadre}\vskip\vcadre}}
\def\arx#1[#2]{\ifcase#1 \relax \or%
  \ar @{-}[#2]  \or%
  \ar @2{-}[#2] \or%
  \ar @{--}[#2] \or%
  \ar @2{.}[#2] \or%
  \ar @{~}[#2]  \fi}

\usepackage{tikz}
\usetikzlibrary{arrows.meta}

\newtheorem{example}{Example}[section]

\newtheorem{theorem}[example]{Theorem}

\newtheorem{corollary}[example]{Corollary}
\newtheorem{conjecture}[example]{Conjecture}
\newtheorem{question}[example]{Question}

\newtheorem{proposition}[example]{Proposition}

\newtheorem{lemma}[example]{Lemma}

\def\Proof{\noindent \it Proof -- \rm}
\def\qed{\hspace{3.5mm} \hfill \vbox{\hrule height 3pt depth 2 pt width 2mm}
\bigskip}

\def\ev{{\rm ev}}

\def\std{{\rm std}}

\def\P{{\bf P}}

\def\SG{{\mathfrak S}}

\def\Pic{{\rm Peak\ }}
\def\Pin{{\rm Pin\ }}
\def\P{{\mathfrak P}}

\title{Pinnacle sets revisited}

\author[J. Falque, J.-C.~Novelli and J.-Y.~Thibon]%
{Justine Falque, Jean-Christophe Novelli and Jean-Yves Thibon}

\address[] {Laboratoire d'Informatique Gaspard Monge, Universit\'e
Gustave-Eiffel, CNRS, ENPC, ESIEE-Paris \\
5 Boulevard Descartes \\Champs-sur-Marne \\77454 Marne-la-Vall\'ee cedex 2 \\
FRANCE}
\email[Justine Falque]{falque@univ-mlv.fr}
\email[Jean-Christophe Novelli]{novelli@univ-mlv.fr}
\email[Jean-Yves Thibon]{jyt@univ-mlv.fr} 



\date{}

\begin{document}

\begin{abstract}
In 2017, Davis, Nelson, Petersen, and Tenner 
[Discrete Math. 341 (2018),3249--3270] initiated the combinatorics of
pinnacles in permutations. We provide a simple and efficient recursion to
compute $p_n(S)$, the number of permutations of $\SG_n$ with pinnacle set
$S$, and a conjectural closed formula for the related numbers $q_n(S)$.
We determine the lexicographically minimal elements of the orbits of the
modified Foata-Strehl action, prove that these elements form a lower ideal of
the left weak order and characterize and count the maximal elements of this
ideal.
\end{abstract}

\maketitle

\section{Introduction}

Given a permutation $\sigma=\sigma_1\dots\sigma_n$, a well-known and
well-studied statistic is the \emph{peak set} of $\sigma$ defined as
\begin{equation}
\Pic \sigma = \{i\,|\, \sigma_{i-1}<\sigma_i>\sigma_{i+1}\}
  \subseteq [2,n-1].
\end{equation}

More recently, some authors~\cite{DNPT} studied a statistic very close in its
definition but rather different in its behaviour: the so-called
\emph{pinnacle set} defined as
\begin{equation}
\Pin \sigma = \{\sigma_i\,|\,\sigma_{i-1}<\sigma_i>\sigma_{i+1}\}
  \subseteq [3,n].
\end{equation}

A set $S$ is an \emph{admissible pinnacle set} (or pinnacle set, for short) if
there is a permutation $\sigma$ such that $\Pin\sigma=S$. It is easy to show
that the pinnacle sets are exactly the sets $S=\{s_1<s_2<\dots<s_p\}$
satisfying $s_i>2i$~\cite{DNPT}.

One main question left unsolved is how to compute efficiently the value
$p_n(S)$ defined in~\cite{DLMSSS} (denoted there by $p_S(n)$ but we prefer
the subscript to be an integer and not a set to enhance readability) as:
\begin{equation}
p_n(S) = \#\{\sigma\in\SG_n \,|\, \Pin\sigma = S\}.
\end{equation}

We shall provide an efficient inductive formula to compute $p_n(S)$
(Theorem~\ref{th-recpn}). Indeed, its complexity is polynomial in both $n$ and
$|S|$ (Proposition~\ref{prop-complex}). Fang in~\cite{Fang} has another
strategy to very efficiently compute $p_n(S)$.

\medskip
It has been noted in~\cite{DLHHIN} that $p_n(S)$ is divisible by $2^{n-1-|S|}$
and in this reference the authors constructed a set of representatives,
called minimal elements, of the orbits of a group action very close to
the Foata-Strehl action, called the dual Foata-Strehl action.
We shall also use this construction to provide another set of representatives,
this time the lexicographically smallest elements of each orbit
(Subsection~\ref{sec-Ppn}), show that they form a lower ideal of the left
weak order (Theorem~\ref{ppn-ideal}), and characterize the maximal elements
of this ideal (Theorem~\ref{max-ppn}).

In Section~\ref{sec-qn}, we study a weighted sum $q_n(S)$ of $p_n(S)$
presented in~\cite{DNPT} in order to simplify the general formulas. We
conjecture a closed formula for $q_n(S)$, in the form of expressions depending
on the size of $S$. The formula looks rather surprising and even if it is
possible to prove some special cases when $|S|$ is small, we have not been
able to prove it in general.

In Section~\ref{sec-arbs}, we consider a set which is equinumerous to the
pinnacle sets. We prove a formula and conjecture another one, both quite
easily expressed in this context. We expect that these objects may shed
some light on other enumerative questions about pinnacle sets.

Finally, in Section~\ref{sec-ords}, we prove that the componentwise
comparison order on pinnacle sets is compatible to the evaluation of $p_n$
and disprove a result that was (erroneously) stated in~\cite{DNPT}.

\bigskip
{\bf Acknowlegements.} This research has been partially supported by the
CARPLO project of the Agence Nationale de la recherche (ANR-20-CE40-0007).

\section{Pinnacle sets}

Let $\P_n(S)$ be the set corresponding to the numbers $p_n(S)$ defined
earlier:
\begin{equation}
\P_n(S) = \{\sigma\in\SG_n \,|\, \Pin\sigma = S\}.
\end{equation}
As observed in~\cite{DLMSSS}, it is immediate that if we have
$\max(S)\leq n$ then
\begin{equation}
\label{pnplusun}
p_{n+1}(S) = 2\,p_n(S).
\end{equation}
Indeed, starting with an element $\sigma$ of $\P_{n+1}(S)$, the value $n+1$ is
not a pinnacle, so it must be either at the first or at the last position of
$\sigma$. Then, removing it provides an element of $\P_n(S)$. This map is
surjective on $\P_n(S)$ since given any element of $\P_n(S)$, there are two
ways of going backwards (put back $n+1$ either in first or last position),
whence the coefficient $2$.

\section{An efficient recursion for $p_n(S)$}

We shall make use of the following notations.
Let $S=\{s_1,\dots,s_p=n\}$ be a pinnacle set.
Let $t(S)$ (or $t$ for short when there is no ambiguity) be the largest
value smaller than $s_p$ that does not belong to $S$.
Then write $S$ as $S=T\cup \{t\!+\!1,\dots,n\}$, with
$T=\{s\in S\,|\, s<t\}$.

\begin{theorem}
\label{th-recpn}
Let $S=\{s_1,\dots,s_p=n\}$ be a pinnacle set.
Let $t$ and $T$ be defined as before.

We then have
\begin{equation}
\label{monpn}
\begin{split}
p_n(T\cup \{t\!+\!1,\dots,n\})
    &= p_n(T\cup \{t,t\!+\!2,t\!+\!3,\dots,n\}) \\
    &+ 2\, p_{n-1}(T\cup \{t,t\!+\!1,\dots,n\!-\!1\}) \\
    &+ 2(n-t)\, p_{n-2}(T\cup \{t,t\!+\!1,\dots,n\!-\!2\}).
\end{split}
\end{equation}
\end{theorem}

\medskip
{\footnotesize 
For example,
\begin{equation}
p_{12} (\{4,8,12\})
  = p_{12} (\{4,8,11\}) +
    2\,p_{11}(\{4,8,11\}) +
    2\,p_{10}(\{4,8\})
\end{equation}
and
\begin{equation}
p_{12} (\{4,8,11,12\})
  = p_{12} (\{4,8,10,12\}) +
    2\,p_{11}(\{4,8,10,11\}) +
    4\,p_{10}(\{4,8,10\}).
\end{equation}
}

\Proof
Let $S=\{s_1,\dots,s_p=n\}$ be a pinnacle set and consider the possible
positions of $t$ among the elements of $\P_n(S)$.
We shall show that $\P_n(S)$ can be partitioned into three sets, each one
being in correspondence with (multiple copies of) other sets $\P_{n'}(S')$.

Let $\sigma$ be an element of $\P_n(S)$. Since $t$ is not a pinnacle, it is
either on a side of $\sigma$, or next to at least one value greater than
itself. We partition those cases into three disjoint sets of cases: either
$t$ is not next to $t+1$, or $t$ is on a side (and next to $t+1$), or $t$ has
two neighbours and one is $t+1$. This last case splits into two subcases
depending on the value of the second neighbour of $t$.

\medskip
Case 1 : $t$ is not next to $t+1$.
Since $t$ is not a pinnacle value but $t+1$ is, if one exchanges $t$ and
$t+1$, one gets an element of $\P_n(S')$, where $S'=T\cup
\{t,t+2,t+3,\dots,n\}$: the pinnacle value $t+1$ is transformed into the
pinnacle value $t$ and all other values remain (non)pinnacle as in $\sigma$.
Moreover, this is a bijection onto $\P_n(S')$ since applying the same
elementary transposition, one gets back from $\P_n(S')$ to $\P_n(S)$. This is
the first term in the formula of the statement.

\medskip
Otherwise, $t$ and $t+1$ are next to each other.
Case 2: $t$ has as neighbours $t+1$ and another one smaller than itself.
Then remove $t$ from $\sigma$ and decrement all values from $t+1$ up so that
we end up with a permutation. We then get an element of $\P_{n-1}(S'')$,
with $S''=T\cup \{t,t+1,\dots,n-1\}$.
Now starting with an element of $\P_{n-1}(S'')$, there are two ways to get
back to $\P_n(S)$ with an inverse transform: increment all values from $t$ up
and put $t$ either at the left or at the right of $t+1$, hence the coefficient
$2$ in the second term of the statement.

\medskip
Case 3: $t$ has no other neighbour than $t+1$.
It is at the first or last position of $\sigma$.
If one removes both $t$ and $t+1$, one gets (up to renumbering the
large values) an element of $\P_{n-2}(S''')$,
with $S'''=T\cup \{t,t+1,\dots,n-2\}$.
There are again two ways of going backwards: start with an element $p$ of
$\P_{n-2}(S''')$ and call $p'$ the element obtained by incrementing twice all
values from $t$ up. Then glue either $(t,t+1)$ at the beginning of $p'$ or
$(t+1,t)$ at the end of $p'$.
Since all values greater than $t$ are pinnacles, they cannot be at an
extremity of $p'$ so we indeed get in both cases an element of $\P_n(S)$.
This contributes a coefficient $2$ in front of the third term of the
statement.

\medskip
Case 4: $t$ is next to two elements greater than itself.
Then remove both $t$ and $t+1$.
We then get again an element of $\P_{S'''}(n-2)$ up to shifting values.
We then have to put back $t$ and $t+1$ and since there are $(n-t-1)$ values at
least equal to $t+2$ and two options to put back the pair $(t,t+1)$ next to
any of those (either $t+1,t$ to their left, or $t,t+1$ to their right), this
case contributes a coefficient $2(n-t-1)$ in front of the third term of the
statement.

\medskip
Summing everything up shows the desired formula.
\qed

\subsection{Known formulas for $p_n(S)$}

Let us now recover the formulas when $S$ has at most two elements.

Consider $S=\{n\}$. Then $T=\emptyset$ and $t=n-1$, so that
\begin{equation}
\label{casun}
p_n(\{n\}) = p_n(\{n-1\}) + 2\,p_{n-1}(\{n-1\}) + 2\,p_{n-2}(\{\}),
\end{equation}
which simplifies into
\begin{equation}
p_n(\{n\}) = 4\, p_{n-1}(\{n-1\}) + 2^{n-2},
\end{equation}
and we recover  Formula~(6) of Prop. 3.6 of~\cite{DNPT}:
\begin{equation}
p_n(\{n\}) = 2^{n-2} (2^{n-2}-1).
\end{equation}

\medskip
With $S=\{n-1,n\}$, we have $T=\emptyset$ and $t=n-2$, and 
\begin{equation}
\label{casdeux}
p_n(\{n-1,n\}) = p_n(\{n-2,n\}) + 2\,p_{n-1}(\{n-2,n-1\})
               + 4\,p_{n-2}(\{n-2\}),
\end{equation}
and with $S=\{\ell,n\}$ and $\ell<n-1$, one obtains an equation similar
to~\eqref{casun}:
\begin{equation}
p_n(\{\ell,n\})
  = p_n(\{\ell,n-1\}) + 2\,p_{n-1}(\{\ell,n-1\}) + 2\,p_{n-2}(\{\ell\}).
\end{equation}
Putting the last two equations together, we recover Formula~(7) of~\cite{DNPT}:
\begin{equation}
p_n(\{\ell,n\}) = 2^{2n-\ell-5}\left(3^{\ell-1}-2^\ell+1\right)
   - 2^{n-3}\left(2^{\ell-2}-1\right).
\end{equation}

\subsection{Evaluating $p_n(S)$ with our formula}

Note that the sum in Formula~\eqref{monpn} only contains positive coefficients
and in particular is not obtained by inclusion-exclusion.

Moreover, this formula can be applied recursively in order to compute
$p_n$: the sum of the elements in each set on the right-hand side 
is strictly smaller than the sum of $S$,
hence showing that the iteration of the formula always stops either on a
non-pinnacle set or on the empty set.
In particular, the induction relies only upon the initial values
$p_n(\emptyset)=2^{n-1}$, all derived from $p_1(\emptyset)=1$.

We shall now discuss the complexity of this algorithm. It is surprisingly low.
Define $C(S)$ as the set of all pinnacle sets that will be involved at some
point when iterating Formula~\eqref{monpn} on a set $S$.
The size of $C(S)$ measures the complexity of the computation.

Note that in Formula~\eqref{monpn} the first and second terms have $|S|$
elements in their pinnacle sets whereas the third one has only $|S|-1$.
Note on the other hand that the second term appears at a later inductive step
when computing the first one. So the number of sets of size $|S|$ in $C(S)$
is equal to the number of sets appearing in the simpler recursion
\begin{equation}
p_n(T\cup \{t\!+\!1,\dots,n\}) = p_n(T\cup \{ t,t\!+\!2,t\!+\!3,\dots,n\}).
\end{equation}

But this equation is easy to analyse in terms on complexity: it requires only
one new element at each step, whose sum of values is strictly smaller than the
previous one: it necessarily stops after $\sum s_i$ steps, so as for the
particular sets of $C(s)$ of size $|S|$, their number is smaller than $n|S|$.

Now, regarding the other sets, a quick look at Formula~\eqref{monpn}
shows that the sets of cardinality $|S|-1$ are all called from the set
$p_{n-2}(T\cup \{t,t\!+\!1,\dots,n\!-\!2\})$. This same idea applies
again to sets of smaller cardinality, hence showing that $A(S)$ is at most
$n \binom{|S|}{2}$.

\begin{proposition}
\label{prop-complex}
The total number of pinnacle sets required to compute $p_n(S)$ recursively
is at most $n|S|^2/2$.
\end{proposition}

The complexity of our computation is therefore a (low degree)
\emph{polynomial} in both $n$ and $|S|$.

\medskip
{\footnotesize
The following Python code is very efficient:
\begin{verbatim}
def is_pinset(S):
    S = sorted(S)
    return all([S[i]>2*(i+1) for i in range(len(S))])

pindic = {tuple([]):0} 

def p(S, n):
    S = tuple(sorted(S))
    if (S,n) in pindic: return pindic[S,n]
    if not is_pinset(S): return 0
    if not S: return 2**(n-1)
    if n<=2: return 0
    if n>S[-1]: return 2**(n - S[-1])*p(S,S[-1])
    m = n
    while m in S: m-=1
    sa = tuple([z for z in S if z>m])
    sb = tuple([z for z in S if z<m])
    X = sb + (sa[0]-1,) + sa[1:]
    Y = sb + tuple([z-1 for z in sa])
    Z = sb + tuple([z-1 for z in sa[:-1]])
    res = p(X, n) + 2*p(Y, n-1) + 2*(n-m)*p(Z, n-2)
    pindic[S,n] = res
    return res  
\end{verbatim}
For example,
\begin{verbatim}
>>> p({5,17,31,42,79,88,97},100)
175144760022244699153358193204473616098046926340653078867873357075537934828864484147200
\end{verbatim}
and the number of recursive calls can be read from the dictionary 
\begin{verbatim}
>>> len(pindic)
753
\end{verbatim}
}

\section{A system of subrepresentatives of classes}

\subsection{Grouping permutations with the same pinnacle set}

First note that both $p_n(\{\})$ and $p_n(\{k\})$ are multiples of
$2^{n-1-|S|}$.
And Formula~\eqref{monpn} shows that it is also the case for any $S$:
the first term has the required power of $2$ by the inductive hypothesis and the
second and third terms have an inductive factor of $2^{n-2-|S|}$, but both
have an extra factor $2$ which yields the required power of $2$.

Let us define $p'_n(S)$ as
\begin{equation}
p'_n(S) := p_n(S)/2^{n-1-|S|}.
\end{equation}
Thanks to Formula~\eqref{pnplusun}, the value of $p'_n(S)$ is independent of
$n$ if $n\geq \max(S)$ so the notation $p'(S)$ makes sense without $n$.

Translating Formula~\eqref{monpn} on $p'$s, we get the following
induction:
\begin{equation} 
\label{monppn}
\begin{split}
p'_n(T\cup \{k\!+\!1,\dots,n\})
    &= p'_n(T\cup \{ k,k\!+\!2,k\!+\!3,\dots,n\}) \\
    &+ p'_{n-1}(T\cup \{k,k\!+\!1,\dots,n\!-\!1\}) \\
    &+ (n-k)\, p'_{n-2}(T\cup \{k,k\!+\!1,\dots,n\!-\!2\}).
\end{split}
\end{equation}

This divisibility property was already noted in~\cite{DLHHIN}
where the authors construct a set of representatives  called  minimal elements
of the orbits of their dual Foata-Strehl action (see their Definition 3.5).

\subsubsection{The (modified) Foata-Strehl action of~\cite{DLHHIN}}

Let $\sigma$ be a permutation and $k$ be a letter of $\sigma$.
Define the \emph{factorization of $\sigma$ according to $k$} as
\begin{equation}
\sigma = \alpha_k\, \beta_k\, k\, \gamma_k\, \delta_k,
\end{equation}
where $\beta_k$ (respectively $\gamma_k$) is the longest (consecutive)
sequence of letters smaller than $k$ immediately to the left (resp. right) of
$k$.
When there is no ambiguity, we shall write $\alpha$ instead of $\alpha_k$ and
similarly for the other factors.

If one represents $\sigma$ as a path with successive heights equal to the
values of $\sigma$, the sets $\beta$ and $\gamma$ represent what $k$ ``sees''
below itself, values higher than $k$ ``blocking'' its view.

Then define as in~\cite{DLHHIN} the map $\varphi_k(\sigma)$ as
\begin{equation}
\varphi_k(\alpha_k\, \beta_k\, k\, \gamma_k\, \delta_k)
= \alpha_k\, \gamma_k\, k\, \beta_k\, \delta_k.
\end{equation}

This collection of maps satisfy some simple properties, all direct corollaries
of the original paper of Foata and Strehl~\cite{FS}:
\begin{proposition}
\label{prop-FS}
Let $\sigma$ be a permutation.

\begin{itemize}
\item $\varphi$ does not change the pinnacle set of $\sigma$,
\item $\varphi$ does not change the vale set of $\sigma$, 
the set of values smaller than (both) their neighbour(s).
\item if one defines $\sigma'=\varphi_{\ell}(\sigma)$ and, if one writes
$\sigma=\alpha \beta k \gamma \delta$
and $\sigma'=\alpha' \beta' k \gamma' \delta'$
as their factorizations according to $k\not=\ell$, then $\beta$ and
$\beta'$ have the same values (maybe not in the same order) and the same is
true for $\gamma$ and $\gamma'$.
\item all the $\varphi$s commute so that if $X=\{x_1,\dots,x_j\}$ is a set,
$\varphi_X=\varphi_{x_1}\circ\dots\circ\varphi_{x_j}$ is independent of the
order of the $x$s, hence well-defined,
\item the orbit of a permutation $\sigma$ is of cardinality $2^{n-1-|S|}$
where $S$ is the pinnacle set of $\sigma$.
\end{itemize}
\end{proposition}

The last property is best understood in terms of the vales of $\sigma$.
Indeed, one can check that $\varphi_x(\sigma)=\sigma$ iff $x$
is a vale of $\sigma$ and since the number of vales $|V|$ is the number of
peaks plus one, the cardinality of the orbit is $2^{n-|V|}$.

Then, given a permutation $\sigma$, the authors of~\cite{DLHHIN} define a
particular subset $W(\sigma)$ of values of $\sigma$ and show that any
element $\tau$ in the same orbit as $\sigma$ satisfies that
$e=\varphi_{W(\tau)}(\tau)$ does not depend on $\tau$, and they call
$e$ the FS-minimal element of the orbit.
Note that this definition might be misleading since their elements are not the
lexicographically minimal elements of each orbit.


Since the orbits of the modified Foata-Strehl action have cardinality
$2^{n-1-|S|}$, they split the sets $\P_n(S)$ into $p'_n(S)$ classes. Given
that the induction on $p'$ of Equation~\eqref{monppn} comes from the
induction on $p$ given in~\eqref{monpn} which is itself an induction on
sets, we translate~\eqref{monppn} into an induction on sets $\P'_n(S)$.
We will show later that these sets happen to be a section of the orbits,
and more precisely, the lexicographically minimal elements of each orbit.

\subsection{The set $\P'_n(S)$ enumerated by $p'_n(S)$}
\label{sec-Ppn}

Let $S$ be an admissible pinnacle set and write as before
$S=T\cup \{t+1,t\!+\!2,t\!+\!3,\dots,n\}$ where $t$ is $t(S)$.

First, as an initialization of the induction, define $\P'_n(\emptyset)$ as
the identity permutation $1\dots n$.
More generally, the set $\P'_s(S)$ with $s>n$ is obtained from
$\P'_{s-1}(S)$ by adding $s$ to the right of all its elements.

Finally, if $n=\max(S)$, the set $\P'_n(S)$
is the union of the following three sets:
\begin{itemize}
\item First, consider the set $A_1=\P'_n(T\cup \{t,t\!+\!2,t\!+\!3,\dots,n\})$.
Then our first subset of $\P'_n(S)$ is obtained by exchanging $t$ and $t+1$ in
each element of $A_1$.
\item Second,
consider the set $A_2=\P'_{n-1}(T\cup \{t,t\!+\!1,\dots,n\!-\!1\})$.
Then our second subset of $\P'_n(S)$ is obtained by sending each element of
$A_2$ to that obtained by mapping its values $x$ onto
\begin{equation}
\left\{
\begin{tabular}{lr}
$x$        & \text{\ if\ } $x<t$, \\
$t,t+1$\   & \text{\ if\ } $x=t$, \\
$x+1$      & \text{\ otherwise.}
\end{tabular}
\right.
\end{equation}
\item Third,
consider the set $A_3=\P'_{n-2}(T\cup \{t,t\!+\!1,\dots,n\!-\!2\})$.
Then our third subset of $\P'_n(S)$ is the union of $n-t$ different sets
obtained from $n-t$ different maps from $A_3$.

For each value $q$ between $t+1$ and $n-1$, send an element of $A_3$ to that
obtained by mapping its values $x$ onto
\begin{equation}
\left\{
\begin{tabular}{lr}
$x$          & \text{\ if\ } $x<t$, \\
$t+1,t,x+2$\ & \text{\ if\ } $x=q$, \\
$x+2$        & \text{\ otherwise.}
\end{tabular}
\right.
\end{equation}

The last set is obtained by mapping each element of $A_3$ to that 
obtained by mapping its values $x$ onto
\begin{equation}
\left\{
\begin{tabular}{lr}
$x$        & \text{\ if\ } $x<t$, \\
$x+2$      & \text{\ otherwise}
\end{tabular}
\right.
\end{equation}
and gluing $t+1,t$ to the right of it.
\end{itemize}

\medskip
{\footnotesize
For example, one can compute that $\P'_6(\{4,6\})$ is
\begin{equation} 
\{ 124365, 134265, 142365, 142563, 142635, 143562, 143625, 156243, 162435 \},
\end{equation}
that $\P'_5(\{4,5\})$ is
\begin{equation} 
\{14253, 14352, 15243\},
\end{equation}
and that $\P'_4(\{4\})$ is
\begin{equation} 
\{ 1243, 1342, 1423\},
\end{equation} 
so that we get the union of the following three sets as
$\P'_6(\{5,6\})$ (in that case $t=4$):
\begin{equation} 
\begin{split}
  A_1 &:= \{125364, 135264, 152364, 152463, 152634, 153462, 153624, 146253,
            162534 \}, \\
  A_2 &:= \{145263, 145362, 162453 \}, \\
  A_3 &:= \{125463, 135462, 154623 \} \cup \{126354,136254,162354\}. \\
\end{split}
\end{equation}
}

\subsection{Properties of $\P'_n(S)$}

\subsubsection{Cardinality of $\P'_n(S)$}

The first property of $\P'_n(S)$ is that it is indeed a subset of the
set of permutations of $\SG_n$ with pinnacle set $S$.

\begin{lemma}
Given a set $S$ and an integer $n$, the elements of $\P'_n(S)$ all have $S$ as
pinnacle set and are distinct. Moreover, $\#\P'_n(S)=p'_n(S)$.
\end{lemma}

\Proof
The facts that the pinnacle set of any element of $\P'_n(S)$ is $S$ and that
there are no repetitions in $\P'_n(S)$ are immediate by induction thanks to
the argument of the proof of Theorem~\ref{th-recpn}: given any element of
$\P'_n(S)$, one can easily find to which case it corresponds, and compute the
representative of $\P'_{n'}(S')$ it comes from. Note also that all maps
defining the third subset of $\P'_n(S)$ have disjoint images. All maps being
disjoint, the cardinalities of $\P'$ follow from the induction on
$p'$s.
\qed

\subsubsection{Characterization of the elements $\P'_n(S)$}


\begin{proposition}
\label{prop-caracR}
A permutation $\sigma\in\SG_n$ with pinnacle set $S$ belongs to $\P'_n(S)$ iff
\begin{itemize}
\item it begins with $1$,
\item it has no double descents (no position $i$ such that
$\sigma_{i-1}>\sigma_i>\sigma_{i+1}$),
\item if $k$ is a pinnacle, $\min(\beta_k)<\min(\gamma_k)$.
\end{itemize}
\end{proposition}

In our representation of permutations as landscapes, this last condition
translates as: the smallest vale that a pinnacle ``sees'' is to its left.

\Proof
Our proof decomposes into three steps. First, we show that all elements of
$\P'_n(S)$ satisfy the properties of the statement. Second, we show that if
$\sigma\in\P_n(S)$ and satisfies the conditions of the statement
then its preimage $\sigma'$ in the inductive definition of $\P_n(S)$ also
does and hence belongs to a $\P'_{n'}(S')$ by induction.
Third, we show that starting from $\sigma'$ back to $\P_n(S)$, the only way
to get an element satisfying the conditions is to use the maps defining
$\P'_n(S)$.

\medskip
We shall first see that all elements of $\P'_n(S)$ satisfy the properties of
the statement.
Indeed, the elements of $\P'_n(S)$ all begin with $1$ and have no double
descents directly by their inductive definition. Regarding the property on
peaks, we do not add peaks in the first or second case of the induction so the
property holds. Let us now consider an element $r$ of $A_3$. If one modifies
some of its values and adds $t+1,t$ to its right, the property holds for all
pinnacles of $r$: either they still view the same elements to their left
and right, or they now see the extra elements $t$ and $t+1$ to their right.
The new pinnacle $t+1$ also satisfies the condition since all vales are
smaller than $t$, it sees one vale to its left (which is smaller than
$t$) and it only sees $t$ to its right.

Now, if one glues $(t+1,t)$ next to a pinnacle greater than $t+1$ in $r$,
the property still holds for the same reasons as before: either a pinnacle
sees the same elements to its left and right as in $r$ or it sees $t$ and
$t+1$ as extra elements one side or the other. But since it also sees other
vales that are by definition of $t$ smaller than $t$, seeing $t$ and $t+1$ is
irrelevant to knowing which side the minimum is.

\medskip
Conversely, let us show that any element $\sigma$ satisfying the
conditions belongs to a $\P'_n(S)$ by induction on its pinnacle set $S$.
Define $t(S)$ as usual on pinnacle sets and find out which case applies to
$\sigma$. Then apply the transform required to get its preimage $\sigma'$ (see
the proof of Theorem~\ref{th-recpn}).
Then $\sigma'$ satisfies the conditions.
Indeed, in Case 1, exchanging $t$ and $t+1$ if they are not neighbours does
not change the position of $1$, does not create a double descent, and does not
change the relative order on minima around pinnacles since $t+1$ was not
a vale to begin with (and $t$ and $t+1$ are consecutive).
In Case 2, removing $t$ if it is next to $t+1$ on one side and to a value
smaller that itself on the other again does not move $1$, does not create a
double descent, and since $t$ was not a vale, the minima change nowhere. In
Cases 3 and 4, removing $t$ and $t+1$ if they are next to a higher pinnacle
does not move $1$, does not create a double descent, and all minima stay the
same. So in all cases, the pre-image element $\sigma'$ also satisfies the
conditions and hence belongs to a $\P'$ by induction.

\medskip
Let us now prove that going back from $\sigma'$ to an element of $\P_n(S)$ 
satisfying the conditions of the statement can only be done using the maps
in the definition of $\P'_n(S)$.
Indeed, let us use the notations of this definition and assume
that $\sigma'$ is in $A_1$. There is only one way to get an element of
$\P_n(S)$ so there is nothing to prove.
Now, if $\sigma'$ is in $A_2$, one gets back to $\P_n(S)$ by adding $t$ either
to the left or to the right of $t+1$ (Case $2$ in the proof of
Theorem~\ref{th-recpn}). But gluing $t$ to the right of $t+1$ does not work
since this would create a double descent.
Finally, if $\sigma'$ is in $A_3$, we have to put back $t$ and $t+1$ together,
either at an extremity of $\sigma'$ or next to a higher pinnacle.
Putting $(t+1,t)$ at the beginning of $\sigma'$ does not work since in that
case the pinnacle $t+1$ would violate the third condition. And putting
$(t,t+1)$ to the right of a higher pinnacle would also not work for the same
reason. So the only inductive steps moving back to an element of $\P_n(S)$
satisfying the conditions of the statement require to start from an element of
a $\P'$ and only apply the induction steps defining $\P'$, hence the result.
\qed

\subsubsection{$\P'_n(S)$ and the orbits of FS} 

We shall make use of the lexicographic order on permutations and denote it by
$\le$.

\begin{lemma}
\label{lem-descorbite}
Let us consider an element $\sigma$ in an orbit of the modified Foata-Strehl
action.

If $\sigma$ has a double descent $\sigma_{i-1}>\sigma_{i}>\sigma_{i+1}$, then
$\varphi_{\sigma_i}(\sigma)<\sigma$.

If $\sigma$ has no double descent and if there are pinnacles $k$
in $\sigma$ so that $\min(\beta_k)>\min(\gamma_k)$ in their factorization, let
$\ell$ be the smallest such pinnacle. Then $\varphi_\ell(\sigma)<\sigma$.
\end{lemma}

\Proof
The case of the double descent is immediate since $\beta$ is empty and all
letters in $\gamma$ are smaller than $\sigma_i$ and move to its left.

Let $\ell$ be defined as in the statement. Then starting from it and moving
left, it sees a nonempty succession of vales with pinnacle smaller than itself
in between until it meets a higher pinnacle.
Since these intermediate pinnacles $\ell'$ are smaller than $\ell$, they
satisfy $\min(\beta_{\ell'})<\min(\gamma_{\ell'})$ so that the vales are
decreasing. Moreover, $\beta_\ell$ cannot begin with a descent since either
it begins with $1$ or it has a pinnacle before it and we assumed that $p$ had
no double descents. So the first letter of $\beta_\ell$ is its minimum.
The same holds to the right of $\ell$ in $p$, so the first letter of
$\gamma_\ell$ is its minimum. Therefore $\varphi_\ell(p)<p$.
\qed

\begin{theorem}
The elements of the sets $\P'_n(S)$ are the lexicographically minimal elements
of their orbits.
\end{theorem}

\Proof
Thanks to Lemma~\ref{lem-descorbite}, we know that an element that does
not satisfy the conditions of~\ref{prop-caracR} cannot be lexicographically
minimal in its orbit since one can apply $\varphi$ to it and obtain a smaller
element. So the lexicographically minimal element of its orbit satisfies the
conditions and thus belongs to a $\P'(S)$.

Now let $\sigma$ be the lexicographically minimal element of an orbit and
consider any other element $\sigma'$ of this orbit.
Since the orbit is connected, there is a set $X$ such that $\varphi_X(p)=p'$.
We can assume that $X$ is minimal, so that any element of $X$ acts non
trivially on $\sigma$. Compute $\sigma''=\varphi_{x_1}(\sigma)$ with
$x_1=\min(X)$.
Then $\sigma''$ violates the conditions next to $x_1$: this $x_1$ cannot be a
vale since $\sigma''\not=\sigma$ so either $x_1$ was in the middle of a double
rise in $\sigma$ and it is a double descent in $\sigma''$ or $\varphi_{x_1}$
exchanged the non-trivial sets $\beta_{x_1}$ and $\gamma_{x_1}$ and hence
violates the third condition.
Now, with all the remaining steps from $\sigma''$ to $\sigma'$, the
neighbourhood of $x_1$ does not change (see Proposition~\ref{prop-FS}) so
$\sigma'$ also violates the conditions on $x_1$.
So no other element in the orbit satisfies all conditions but
the lexicographically minimal one.
\qed

Note that our algorithm provides an efficient way to build all
lexicographically minimal elements of the dual Foata-Strehl orbits.

\subsubsection{$\P'_n$ and the left weak order}

Let us consider the whole set $\P'_n$ which is the union of all $\P'_n(S)$
with $S$ a pinnacle set with maximum at most $n$.

Recall that the left weak order on $\SG_n$ is the transitive closure of the
relation on permutations given by $\sigma<\tau$ if $s_i.\sigma=\tau$ and
$\ell(\tau)=\ell(\sigma)+1$ where $s_i$ is the transposition $(i,i+1)$ and
$\ell$ is the number of inversions of permutations.
And recall that a lower ideal of a poset is a subset $S$ of this poset such
that $s\in S$ implies that $t\in S$ for all $t<s$.

\begin{theorem}
\label{ppn-ideal}
$\P'_n$ is a lower ideal of the left weak order.
\end{theorem}

\Proof
We just have to prove that if $\tau\in\P'_n$, then $s_i.\tau=\sigma\in\P'_n$
too if $\ell(\tau)=\ell(\sigma)+1$.
Thanks to their characterization, we know that $\tau$ begins with $1$, has no
double descent and that $\beta_k<\gamma_k$ for any $k$ in the pinnacles of
$\tau$.

The first two criteria are automatic with $\sigma$. Concerning the minima,
since the exchange of $i$ and $i+1$ moves $i+1$ further right and since
$i$ and $i+1$ are consecutive values, there is no situation where a left
minima goes from smaller to greater than a right minima.
\qed

When given an ideal, it is customary to consider its maximal elements.
In our case, the first maximal elements of our ideal are
\begin{equation}
\label{minmax}
\begin{split}
& 12 \\
& 132 \\
& 1342, 1423 \\
& 14352, 14523, 15243 \\
& 145362, 153462, 154623, 156243, 162453, 162534 \\
& 1546372, 1563472, 1564723, 1635472, 1645723, \\
& 1657243, 1672453, 1672534, 1725463, 1725634, 1726354 \\
\end{split}
\end{equation}

Recall that the standardisation process $\std$ of a word without repetition
amounts to renumbering the values with $1$ up to its size in the order they
were in the beginning. For example, $\std(15726)=13524$.

\begin{theorem}
\label{max-ppn}
The maximal elements of $\P'_n$ are the elements $\sigma$
of $\P'_n$ satisfying the extra conditions:
\begin{itemize}
\item Cut $\sigma$ as $\sigma=u.n.v$. Then $\std(u)$ and $\std(v)$ are
themselves maximal elements,
\item the values of $v$ form an interval $[2,\ell]$.
\end{itemize}
\end{theorem}

\Proof
Let us first prove that an element $\sigma$ of $\P'_n$ that does not satisfy
one extra condition cannot be maximal.
First, if there is a letter $i>1$ in $u$ such that $i+1$ is in $v$, then
exchange $i$ and $i+1$. This element is still in $\P'_n$ since the only
condition that could fail is the pinnacle condition on $n$, but
$n$ still sees $1$ on the left.
Now, if $\sigma$ satisfies the first extra condition but not the second,
$v$ is composed of consecutive letters and if its standardized
in not maximal, the same transposition (shifted by one) applied to it shows
that $v$ is not maximal either. The same argument applies to $u$.

Conversely, let us prove that an element that satisfies both extra
conditions is indeed maximal. Since $u$ and $v$ are maximal, there is no
transposition inside $u$ or $v$ that could bring another element in $\P'_n$
since violating a condition is independent from what happens on the other side
of $n$. It is also impossible to move $1$ so the only allowed
transposition increasing the inversion number of $\sigma$ is the transposition
$(n-1,n)$. But this one fails since after the exchange $n-1$ sees $2$ on its
right and does not see $1$ anymore on its left.
\qed

\begin{corollary}
Let $M_n$ be the set of permutations defined inductively by
\begin{itemize} 
\item $M_1$ is $\{1\}$,
\item $M_2$ is $\{1,2\}$,
\item $M_n$ is obtained as the union for all $k$ of the sets
$M'_{n,k}$ where an element of $M'_{n,k}$ is the concatenation
$a_{n-1-k} n b_{k}$ where $a_{n-1-k}$ is an element of $M_{n-1-k}$ whose
values at least $2$ have been shifted by $k$ and $b_{k}$ is an element of
$M_k$ whose values have been shifted by $1$.
\end{itemize}

Then $M_n$ is the set of the maximal elements of the ideal of the minimal
elements of the (modified) Foata-Strehl orbits.
\end{corollary}

\begin{corollary}
The maximal elements of $\P'_n$ are enumerated by Sequence~A007477
of~\cite{Slo} (up to a change of indices) since they satisfy the
induction formula
\begin{equation}
p'_n = \sum_{k=1}^{n-2} p'_k p'_{n-1-k}.
\end{equation}
\end{corollary}

Here are the first terms of A007477.
\begin{equation}
1,1,1,2,3,6,11,22,44,90,187
\end{equation}
and one can check that this is consistent with the list given
in~\eqref{minmax}.

\Proof
Immediate from the previous characterization by summing over the different
positions that $n$ can occupy.
\qed

For example, one easily gets $M_4=\{1342,1423\}$ and
$M_5=\{14352,14523,15243\}$ so that $M'_{10,4}$ is obtained as all
concatenations of an element of $\{18796,18967,19687\}$ with $10$ and
with an element of $\{2453,2534\}$.

\section{A conjectural formula for $q_n(S)$}
\label{sec-qn}

In~\cite{DNPT}, Question 4.5, it was conjectured that
\begin{equation}
q_n(S) := \sum_{I\subset S} 2^{|I|} p_n(I)
\end{equation}
had a nice formula and indeed, we shall conjecture a general formula for
$q_n(S)$.

We shall write it as a product $q'_n(S)r_n(S)$,
where $q'$ takes into account (almost all) the powers of $2$ that appear as
factors in the overall formula, and $r$ has a more complicated formula.

Let $S=\{s_1,\dots,s_p\}$ be a pinnacle set and $n\geq s_p$.
Then let
\begin{equation}
q'_n(S) = 2^{n-1} 2^{s_p-2} 2^{-s_{p-1}} 2^{s_{p-2}-2} 2^{-s_{p-3}} \dots
\end{equation}

We shall describe a simple but quite surprising algorithm to compute the
second factor $r_n(S)$.
First, define
\begin{equation}
D(S)=(d_1,\dots,d_{p-1}) := (s_1-1,s_2-s_1,\dots,s_{p-1}-s_{p-2}).
\end{equation}

Then build the following abstract expressions
\begin{equation}
E_1 := x_{1,3} + x_{1,1},
\end{equation}
and
\begin{equation}
E_{2k+1} := f_e(E_{2k}) \text{\ \ and\ \ } E_{2k} := f_o(E_{2k-1}),
\end{equation}
where $f_o$ substitutes all elements involving the largest index $a$ as
\begin{equation}
f_o(x_{a,k}) =
x_{a,k} \left(x_{a+1,\frac{k+1}2} +
              x_{a+1,\frac{k-1}2}
       \right)
\end{equation}
with the convention $x_{a,0}=0$,
and $f_e$ substitutes all elements involving the largest index $a$ as
\begin{equation}
f_e(x_{a,k}) =
x_{a,k} (x_{a+1,2k+1}+x_{a+1,2k-1}).
\end{equation}

{\footnotesize 
For example, we get the following first expressions
\begin{equation}
E_2 = f_o(E_1) = x_{1,3}(x_{2,2}+x_{2,1}) + x_{1,1}x_{2,1},
\end{equation}
\begin{equation}
E_3 = f_e(E_2) = x_{1,3}(x_{2,2}(x_{3,5}+x_{3,3}) + x_{2,1}(x_{3,3}+x_{3,1}))
                 + x_{1,1}x_{2,1}(x_{3,3}+x_{3,1}),
\end{equation}
\begin{equation}
\begin{split}
E_4 = f_o(E_3)
   =&\ x_{1,3}\big(x_{2,2}(x_{3,5}(x_{4,3}+x_{4,2})
                   + x_{3,3}(x_{4,2}+x_{4,1})) \\
    &\ \ \ \ \ \ \ + x_{2,1}(x_{3,3}(x_{4,2}+x_{4,1})+x_{3,1}x_{4,1})\big) \\
    &+ x_{1,1}x_{2,1}(x_{3,3}(x_{4,2}+x_{4,1}) + x_{3,1}x_{4,1}).
\end{split}
\end{equation}
}

Then, if $S$ has $p$ elements so that $D(S)$ has $p-1$, define
\begin{equation}
r_n(S) := \ev(E_{p-1}),
\end{equation}
where $\ev$ evaluates $x_{i,j}$ to $j^{d_{p-i}}$.

\begin{conjecture}
For all pinnacle set $S$ and all $n\geq max(S)$, the value of $q_n(S)$ is
equal to the product $q'_n(S)r_n(S)$.
\end{conjecture}

{\footnotesize
For example, if one replaces the $d$s by their values as parts of $S$,
the first formulas read
\begin{equation}
q_n(\{m\}) = 2^{n-1} 2^{m-2},
\end{equation}
\begin{equation}
q_n(\{\ell,m\}) = 2^{n-1 + m-2 - \ell} \left(3^{\ell-1}+1\right),
\end{equation}
\begin{equation}
q_n(\{k,\ell,m\})
 = 2^{n-1 + m-2 - \ell + k-2}
   \left[3^{\ell-k}(2^{k-1}+1) + 1\right],
\end{equation}
\begin{equation}
q_n(\{j,k,\ell,m\})
 = 2^{n-1 + m-2 - \ell + k-2 - j}
   \left[3^{\ell-k}\left(2^{k-j}(5^{j-1}+3^{j-1}) + 3^{j-1}+1\right)
   + 3^{j-1}+1\right],
\end{equation}
\begin{equation}
\begin{split}
q_n(\{i,j,k,\ell,m\})
 =&\ 2^{n-1 + m-2 - \ell + k-2 - j + i-2} \\
 &\ \ \ \big[3^{\ell-k}
           \left(2^{k-j}(5^{j-i}(3^{i-1}+2^{i-1})+3^{j-i}(2^{i-1}+1))
                 + 3^{j-i}(2^{i-1}+1)+1\right) \\
 &\ \ \ \  + 3^{j-i}(2^{i-1}+1)+1\big].
\end{split}
\end{equation}
}

\subsection{Properties of our conjectural formula for $q_n(S)$}

\subsubsection{The variables $x_a$}

By definition of $f_e$, any $E_{2k+1}$ only contains factors $x_{2k+1,i}$
with odd $i$ so that applying $f_o$ to these will only put integer-valued
$j$ in the new variables $x_{2k+2,j}$.

\subsubsection{Structure of the formula}

Thanks to its structure, it is obvious that $E_k$ is a sum of monomials
without multiplicities. 
Moreover, a simple induction shows that $E_k$ exactly has
$\binom{k+1}{\lfloor{(k+1)/2}\rfloor}$ terms,
so that these terms are in bijection with the pinnacle sets with maximum at
most $k+2$.

There is a simple bijection between both sets where the various $x_{k}$s
encode the size of the pinnacle set but this bijection does not seem very
relevant to better understanding (or proving) the formula.

\subsubsection{Special values}

By definition, $q_n(S)$ can be nonzero even if $S$ is not a
pinnacle set. In particular, the definition implies that
\begin{equation}
q_n(\{2\}\cup S) = q_n(S).
\end{equation}

Even if our formula does not coincide in all non-pinnacle cases, it indeed
does seem to coincide in that case. It was in fact a big help in finding the
first formulas with small values of $|S|$.

\subsubsection{Computing $p_n(S)$}

Assuming our formula for $q_n(S)$ is correct, it provides another algorithm to
compute $p_n(S)$: first, compute the expressions $E_k$ up to $k=|S|-1$, then
apply these to all subsets of $S$ and apply an inclusion-exclusion process
to get $p_n(S)$.

Note that it is very possible to brute-force $p_n$ on pinnacle sets with $3$
or $4$ values, get a general formula, and prove it using our induction.
However, we do not expect this strategy to help prove our conjectural formula
for $q_n$ in general.

In terms of computation time, this algorithm is much less efficient that our
algorithm using Formula~\eqref{monpn} in general since it requires to compute
$2^{|S|}$ terms when our first algorithm requires $n|S|^2$ terms.

But we can still expect that a proof and a better understanding of our formula
for $q_n(S)$ could bring ideas to get a general formula for $p_n(S)$.

\section{Ordered forests of complete binary trees}
\label{sec-arbs}

Among all combinatorial objects enumerated by the central binomial sequence
A001405 of ~\cite{Slo} are the pinnacle sets and ordered forests of
complete binary trees.
The bijection between both is very simple and sheds light on some structures
of pinnacle sets.

Let $S=\{s_1,\dots,s_p\}$ be a pinnacle set with $n\geq s_p$.
We shall then build a sequence of complete binary trees as follows:
\begin{itemize}
\item set $k=n$ and start with a sequence of one tree: the tree with one
node.
\item Put $k$ in the rightmost empty node.
If $k\in S$ then draw both children of $k$ as empty nodes.
Otherwise, do not draw those children and if there is no empty node anymore,
put a new tree with one node to the left of the previous one.
Set $k:=k-1$ and start again this step until $k=0$.
\end{itemize}

\medskip
{\footnotesize
For example, if $S=\{4,8,9,12\}$ and $n=13$, we get the following forest:
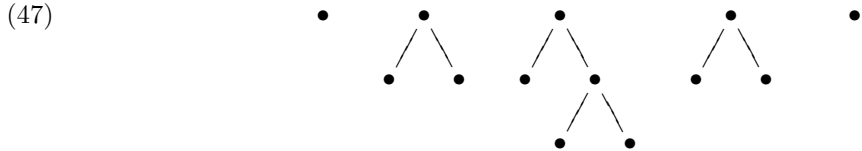
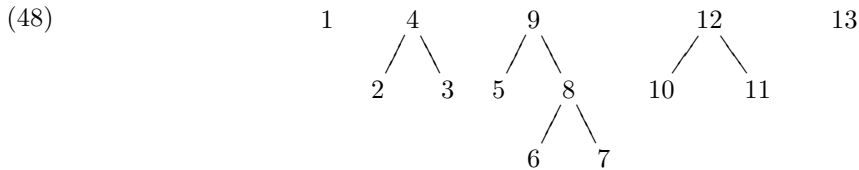
\begin{equation}
\setlength\unitlength{1.7mm}
\xymatrix@R=0.5cm@C=.9mm{
&  *{\GrTeXBox{\bullet}} \\
}
\hskip.5cm
\xymatrix@R=0.5cm@C=.9mm{
&  *{\GrTeXBox{\bullet}}\arx1[ld]\arx1[rd] \\
   *{\GrTeXBox{\bullet}} & &
   *{\GrTeXBox{\bullet}} \\
}
\hskip.5cm
\xymatrix@R=0.5cm@C=.9mm{
&  *{\GrTeXBox{\bullet}}\arx1[ld]\arx1[rd] \\
   *{\GrTeXBox{\bullet}} & &
   *{\GrTeXBox{\bullet}}\arx1[ld]\arx1[rd] \\
&  *{\GrTeXBox{\bullet}} & &
   *{\GrTeXBox{\bullet}} \\
}
\hskip.5cm
\xymatrix@R=0.5cm@C=.9mm{
&  *{\GrTeXBox{\bullet}}\arx1[ld]\arx1[rd] \\
   *{\GrTeXBox{\bullet}} & &
   *{\GrTeXBox{\bullet}} \\
}
\hskip.5cm
\xymatrix@R=0.5cm@C=.9mm{
&  *{\GrTeXBox{\bullet}} \\
}
\end{equation}
or, with labels
\begin{equation}
\setlength\unitlength{1.7mm}
\xymatrix@R=0.5cm@C=.9mm{
&  *{\GrTeXBox{1}} \\
}
\hskip.3cm
\xymatrix@R=0.5cm@C=.9mm{
&  *{\GrTeXBox{4}}\arx1[ld]\arx1[rd] \\
   *{\GrTeXBox{2}} & &
   *{\GrTeXBox{3}} \\
}
\hskip.3cm
\xymatrix@R=0.5cm@C=.9mm{
&  *{\GrTeXBox{9}}\arx1[ld]\arx1[rd] \\
   *{\GrTeXBox{5}} & &
   *{\GrTeXBox{8}}\arx1[ld]\arx1[rd] \\
&  *{\GrTeXBox{6}} & &
   *{\GrTeXBox{7}} \\
}
\hskip.3cm
\xymatrix@R=0.5cm@C=.9mm{
&  *{\GrTeXBox{12}}\arx1[ld]\arx1[rd] \\
   *{\GrTeXBox{10}} & &
   *{\GrTeXBox{11}} \\
}
\hskip.3cm
\xymatrix@R=0.5cm@C=.9mm{
&  *{\GrTeXBox{13}} \\
}
\end{equation}
}

\medskip
Note that by construction the elements of $S$ are exactly the labels of the
internal nodes of the forest.

Conversely, starting with a forest, label the nodes in left suffix order
(order the trees from the left one to the right one, and within a tree, label
first the left subtree of a node, then its right subtree, then the node
itself). Then read the values of the internal nodes.

This is a bijection since each step can easily be reverted between the
pinnacle sets and the left-suffix labelled trees.

The simplest induction on $p_n(S)$ of Equation~\eqref{pnplusun} translates
on these objects as:
\begin{equation}
p_n(F,O) = 2 p_n(F),
\end{equation}
for any forest (sequence of trees) $F$ and where $O$ is the one-node tree.

With the help of this notation, we found two properties on $p_n(S)$.

\begin{proposition}
Let us consider a sequence of trees $F := (T_1,T_2,\dots,T_r)$ encoding a
pinnacle set.

Then
\begin{equation}
p_n(F) = p_n(T_1)\, p_n(O,T_2,\dots,T_r).
\end{equation}
\end{proposition}

\Proof
The values in the tree $T_1$ are necessarily consecutive in any element of
$p_n(S)$. So we go from $p_n(F)$ to $p_n(O,T_2,\dots,T_r)$ by replacing this
sequence by $1$ and standardizing the result. The converse operation changes
$1$ into all possibilities for $p_n(T_1)$ hence explaining the multiplicative
factor.
\qed

\begin{conjecture}
Let us consider a sequence of trees $F := (O,T_2,T_3,\dots,T_r)$ encoding a
pinnacle set.
Then
\begin{equation}
p_n(F) = p_n(O,T_2)\, p_n(O,O,T_3,\dots,T_r).
\end{equation}
\end{conjecture}

A nice proof of this conjecture would probably be a first step before being
able to generalize it to longer sequences of one-node trees at the beginning.

\begin{question}
When the first two trees of $F$ are the tree with one node, $p_n$ does not
factorize. However, we expect there should exist a generalization of this
result as an additive formula.
\end{question}

\section{Orders on pinnacle sets}
\label{sec-ords}

In~\cite{DNPT} Question 4.3 was whether there is a nontrivial order on
sequences of a given size such that $S_1<S_2$ would imply $p_n(S_1)<p_n(S_2)$.

Thanks to an argument we have already seen, it is clear that if $S_1$ is
componentwise smaller than $S_2$, then we indeed have $p_n(S_1)<p_n(S_2)$ (it
seems that the authors had seen that but did not write the proof).
We just need to prove the property if
$S_1=\{s_1,\dots,s_p\}$ and $S_2$ is obtained by changing $s_i$ into $s_{i}+1$
(if $s_{i}+1\not\in S_1$).

Starting with an element of $\P_n(S_1)$, exchange the values $s_i$ and $s_i+1$.
Then we get an element of $\P_n(S_2)$: indeed, since $s_i+1$ was not a
pinnacle (not in $S_1$), it was not next to $s_i$ but was either at an
extremity or next to a greater value. After the exchange, all (non) pinnacles
remain (non) pinnacles except $s_i$ that changes status with $s_i+1$. Note
that the converse is not true (see the explanations of our inductive formula
on $p_n(S)$): in $\P_n(S_2)$, there are elements where $s_i$ is next to
$s_i+1$ and these are not in the image of the previous map.

To get more precise comparison results on pinnacle sets, we would need to
inherit from small cases to larger ones, as in \emph{e.g}, Proposition 3.9
in~\cite{DNPT}.
Unfortunately, the provided proof is incorrect since the induction formula
used to justify it does not hold in general (and would give a direct and easy
formula for $p_n(S)$).
There is no real hope to patch it since
\begin{equation}
p_9(\{3,9\})=1984 > p_9(\{5,6\})=1152
\end{equation}
whereas
\begin{equation}
p_{12}(\{3,9,10,11,12\})=172800 < p_{12}(\{5,6,10,11,12\}) = 207360.
\end{equation}

\section{Concluding remarks}

We did not address here how to put algebraic structures in the picture of
pinnacle sets but the strategy is the same as in the case of top-descent
values of permutations. One cannot build a Hopf algebra, not even a subalgebra
of the algebra on permutations, but a quotient works as in~\cite{HNTT}.
This will be addressed in a forthcoming paper.

\footnotesize
\bibliographystyle{alpha}

\end{document}